\documentclass{article}
\usepackage{maa-monthly}

\theoremstyle{theorem}
\newtheorem{theorem}{Theorem}

\newtheorem{lemma}[theorem]{Lemma}
\newtheorem{corollary}[theorem]{Corollary}

\theoremstyle{definition}

\begin{document}

\title{The convergence of power matrices }
\markright{The convergence of power matrices}
\author{Vyacheslav M. Abramov}  
\thanks{Former affiliation: School of Mathematics, Monash University, Australia}%
\thanks{Current status: aged pensioner}
\thanks{Address: 24 Sagan Drive, Cranbourne North, Melbourne, Victoria-3977, Australia}
\thanks{Email: vabramov126@gmail.com}

\maketitle

\begin{abstract}
For the class of $d\times d$ matrices $B=[b_{i,j}]$ with complex entries satisfying $\sum_{i=1}^{d}|b_{i,j}|=1$, we provide necessary and sufficient conditions for the convergence of power matrices $B^n$ to a nonzero limit matrix. Our result extends the well-known classic one from the theory of Markov chains. The approach of this note is elementary, self-contained and based on standard matrix operations only.
\end{abstract}

%
\section{Introduction.}
The Perron--Frobenius (PF) theorem (established first by Oskar Perron \cite{P} for matrices with positive entries and then extended by Ferdinand Georg Frobenius \cite{F} for matrices with nonnegative entries) is one of the most important results in the theory of real square matrices.
The PF theorem has various applications that include probability theory, theory of dynamical systems, economics, matrix theory etc \cite{HJ, L, S}. One of the widely used applications of the PF theorem that is related to the theory of discrete Markov chains is as follows.
Let $B$ be a real square matrix with nonnegative entries $b_{i,j}$ satisfying $\sum_{j=1}^{d}b_{i,j}=1$, $i=1,2,\ldots,d$ or $\sum_{i=1}^{d}b_{i,j}=1$, $j=1,2,\ldots,d$. If $B$ is aperiodic and irreducible (e.g., \cite[pp. 11--13]{T} for the definitions), then there exists a positive limit of $B^n$ as $n\to\infty$ (see \cite[Ch. 1, Rel. (19)]{T}). The limiting distributions of irreducible Markov chains as an application of the limit of $B^n$ is given in \cite[p. 15]{T}.

By positive limit matrix we mean a matrix with nonnegative entries, some of which are strictly positive. In particular, if all entries of a matrix $B$ are strictly positive, then it is irreducible and aperiodic, and $B^n$ converges to a limit matrix, the entries of which are strictly positive.

A matrix with nonnegative entries satisfying $\sum_{j=1}^{d}b_{i,j}=1$, $i=1,2,\ldots,d$ or $\sum_{i=1}^{d}b_{i,j}=1$, $j=1,2,\ldots,d$, is called stochastic matrix. Stochastic matrices were originally used by Russian mathematician Andrey Markov \cite{M} to describe transition probabilities in certain probability problems with outcomes depending on the experiment conditions and then widely developed covering many areas of application of Markov chains.

Whereas the theory of stochastic matrices with real entries is well-developed and well-known, not much is known about matrices with complex entries.
The PF theory for matrices with complex entries has been developed very recently \cite{Rugh, Rump, NV} and concerned their spectral properties, but not the limits of their powers. Theorems of PF type for real matrices without sign restrictions were obtained in \cite{Rump1}.

The aim of this note is to obtain a natural extension of the aforementioned known result on matrices with real entries for the case of matrices with complex entries.
Specifically, we consider a $d\times d$ matrix $B$ with complex entries $[b_{i,j}]$ assuming that the matrix of absolute values (denoted by $|B|=[|b_{i,j}|]$) is irreducible, aperiodic and satisfying $\sum_{i=1}^{d}|b_{i,j}|=1$, $j=1,2,\ldots, d$, and
provide the necessary and sufficient conditions, under which the sequence $B^n$ converges to a nonzero limit matrix. A sequence of power matrices with complex entries can arise in various applied problems. As well, the results of the present study can be a subject of further analysis of matrices with complex entries.


We use the following global notation. Given a field $\mathbb{F}$ (which will be equal to $\mathbb{R}$ or $\mathbb{C}$ in this note), and integer $d\geq2$, we define $\mathbb{F}^{d\times d}$ to be the set of $d\times d$ matrices with entries in $\mathbb{F}$.

The results obtained in this note can be derived by different ways. We mention the well-known textbook of Felix R. Gantmacher \cite{G} with a number of important results in Chapter XIII. One of them is \cite[Lemma 2, p. 57]{G} (see also \cite{W}) that makes possible to provide similar study.
The advantage of our approach, however, is that our proof is self-contained and based on elementary matrix operations only.

In Section \ref{S2}, we formulate and prove our main result.
The appendix contains auxiliary statements that are used to prove the required theorem.

\section{Main result and its proof.}\label{S2}
Let $B=[b_{i,j}]\in\mathbb{C}^{d\times d}$, the entries of which satisfy the property $\sum_{i=1}^{d}|b_{i,j}|=1$. We prove the following theorem.

\begin{theorem}\label{t1}
Let $B=[b_{i,j}]\in\mathbb{C}^{d\times d}$, $\sum_{i=1}^{d}|b_{i,j}|=1$, and let $|B|$ be aperiodic and irreducible. 
Then the sequence $B^n$ converges to a nonzero limit irreducible matrix if and only if
\[
b_{i,j}=|b_{i,j}|\frac{\xi_i}{\xi_j},
\]
where $\xi_i\in\mathbb{S}^1$, $i=1,2,\ldots, d$.
\end{theorem}

\begin{proof}
To start our proof, we first need to prove the lemma given below.

\begin{lemma}\label{lem1}
Let $A=[a_{i,j}]\in\mathbb{R}^{d\times d}$ and define $\tilde{A}=[\tilde{a}_{i,j}]\in\mathbb{C}^{d\times d}$ as $\tilde{A}=\Xi A\Xi^{-1}$, where $\Xi$ is an invertible diagonal matrix, $\mathrm{diag}(\Xi)=\big[\xi_i\big]_{i=1}^d$ with $\xi_i\in\mathbb{S}^1$. That is,
\begin{equation}\label{10}
\tilde{a}_{i,j}=a_{i,j}\frac{\xi_i}{\xi_j}, \quad i,j=1,2,\ldots,d.
\end{equation}
Then for all $n\geq1$, the entries of $\tilde{A}^n$ and $A^n$ are related as \eqref{10}. Namely,
\begin{equation}\label{11}
\tilde{a}_{i,j}^{(n)}=a_{i,j}^{(n)}\frac{\xi_i}{\xi_j}, \quad i,j=1,2,\ldots,d.
\end{equation}
where $A^n=[a_{i,j}^{(n)}]$ and $\tilde{A}^n=[\tilde{a}_{i,j}^{(n)}]$.
\end{lemma}

\begin{proof}
Indeed, from $\tilde{A}=\Xi A\Xi^{-1}$ we have $\tilde{A}^n=\Xi A^n\Xi^{-1}$ that leads to \eqref{11}. The lemma is proved.
\end{proof}

\subsection{Proof of the sufficient condition.}
Using the polar system, write $b_{i,j}=|b_{i,j}|\mathrm{e}^{\boldsymbol{i}\theta_{i,j}}$, $i,j=1, 2,\ldots, d$, $\boldsymbol{i}=\sqrt{-1}$. According to the assumption,
\[
b_{i,j}=|b_{i,j}|\frac{\xi_i}{\xi_j}\quad i,j=1, 2,\ldots, d.
\]
where $\xi_i=\mathrm{e}^{\boldsymbol{i}\theta_i}$ and $\theta_{i,j}:=\theta_{i}-\theta_j$, $i,j=1, 2,\ldots, d$.

According to the proof of Lemma \ref{lem1}, $B^n=\Xi|B|^n\Xi^{-1}$, and due to \cite[Ch. 1, Rel. (19)]{T}, if $|B|$ is aperiodic and irreducible, then the sequence $|B|^n$ converges to a nonzero limit. Hence, the sequence $B^n$ converges to a nonzero limit too.

\subsection{Proof of the necessary condition.}
Using the notation
\[
X^{(n)} = \big[x^{(n)}_{i,j}\big]_{i,j=1}^d:=B^n,
\]
note that $X^{(n+1)} = X^{(n)}B$, that is,

\begin{equation}\label{0}
x_{i,k}^{(n+1)}=\sum_{j=1}^{d}x_{i,j}^{(n)}b_{j,k}, \quad i,k=1, 2,\ldots, d.
\end{equation}
By hypothesis, there exists a nonzero limit irreducible matrix $X^{(*)}=\lim_{n\to\infty}X^{(n)}$, the entries of which are denoted by $x_{i,j}^{(*)}$.  Then we obtain $X^{(*)}=X^{(*)}B$, and the system of the equations for the entries is
\begin{equation}\label{0.5}
x_{i,k}^{(*)}=\sum_{j=1}^d x_{i,j}^{(*)}b_{j,k}, \quad i,k=1, 2,\ldots, d.
\end{equation}
It is clear that if \eqref{0.5} is satisfied, then we must also have
\begin{equation*}\label{0.6}
x_{i,k}^{(*)}=\sum_{j=1}^d x_{i,j}^{(*)}x_{j,k}^{(2)}, \quad i,k=1, 2,\ldots, d.
\end{equation*}

Taking absolute values, we have:
\begin{equation}\label{1}
|x_{i,k}^{(*)}|\leq\sum_{j=1}^{d}|x_{i,j}^{(*)}x_{j,k}^{(2)}|, \quad i,k=1, 2,\ldots, d.
\end{equation}
If at least one of the $d^2$ inequalities
\[
|x_{i,k}^{(2)}|\leq\sum_{j=1}^{d}|b_{i,j}b_{j,k}|, \quad i,k=1, 2,\ldots, d
\]
is strict, then $\sum_{i=1}^{d}|x_{i,j}^{(2)}|<1$
for at least one of the indices $j$. Note that given that our hypothesis is true, both $B^2$ and $|B^2|$ must be irreducible too. Then
according to Lemma \ref{lemA0} (applied to $|B^2|$), we have
\[
\sum_{i=1}^{d}|x_{i,j}^{(4)}|<1, \quad j=1, 2,\ldots, d.
\]
By Corollary \ref{lemA}, $(|B^2|)^n$ converges to the zero matrix as $n\to\infty$. Since $|B^{2n}|\leq (|B^2|)^n$ entrywise and $|B^{2n+1}|\leq |B^{2n}||B|$ entrywise, then both $|B^{2n}|$ and $|B^{2n+1}|$ converge to the zero matrix too as $n\to\infty$, and hence all $|x_{i,k}^{(*)}|$ in \eqref{1} are zeros. 
Consequently, $X^{(*)}$
is the zero matrix. This contradicts the hypotheses, so
\begin{equation}\label{1.2}
\big|\sum_{j=1}^{d}b_{i,j}b_{j,k}\big|=|x_{i,k}^{(2)}|=\sum_{j=1}^{d}|b_{i,j}b_{j,k}| \quad \forall i, k\in\{1, 2,\ldots, d\}.
\end{equation}
From \eqref{1.2} we conclude that $b_{i,j}b_{j,k}$, $j=1,2,\ldots, d$, have equal arguments, for all $i, k =1, 2,\ldots, d$. Solving for the arguments, we see that (regardless of the entries being $0$ or
not) there exist phases 
$\xi_1, \xi_2,\ldots, \xi_d\in\mathbb{S}^1$ and $\omega_1, \omega_2,\ldots, \omega_d\in\mathbb{S}^1$
such that $b_{i,j}=|b_{i,j}|\xi_i\omega_j$, and
\begin{equation}\label{1.1}
\begin{aligned}
B=\left(\begin{matrix}b_{1,1} &b_{1,2} &\cdots &b_{1,d}\\ b_{2,1} &b_{2,2} &\cdots &b_{2,d}\\ \vdots &\vdots &\cdots &\vdots\\ b_{d,1} &b_{d,2} &\cdots &b_{d,d}\end{matrix}\right)
=\left(\begin{matrix}\xi_1 & & &\\ &\xi_2 & &\\ & &\ddots &\\ & & &\xi_d\end{matrix}\right)|B|\left(\begin{matrix}\omega_1 & & &\\ &\omega_2 & &\\ & &\ddots &\\ & & &\omega_d\end{matrix}\right).
\end{aligned}
\end{equation}
Taking powers, for the limit one must have $\xi_i\omega_i=1$, $i=1, 2,\ldots, d$, which means that the diagonal matrices are simply unitary.

Hence, for the nonzero matrix-solution $X^{(*)}$, we have the equalities
\begin{equation*}\label{1.3}
x_{i,j}^{(*)}b_{j,i}=|x_{i,j}^{(*)}b_{j,i}|, \quad i,j=1, 2,\ldots, d,
\end{equation*}
that in turn follow from
\[
x_{i,j}^{(n)}b_{j,i}=|x_{i,j}^{(n)}b_{j,i}|, \quad i,j=1, 2,\ldots, d, \quad n\geq1
\]
as a consequence of \eqref{1.1} and Lemma \ref{lem1}. \end{proof}

\appendix
\section{Auxiliary statements.}
We prove the following results, where $d$ is a positive integer.
\begin{lemma}\label{lemA0}
Let $B\in\mathbb{R}^{d\times d}$ be an irreducible matrix with nonnegative real entries $b_{i,j}$ satisfying the system of inequalities $\sum_{i=1}^{d}b_{i,j}\leq1$. If at least one of these inequalities is strict, then for the entries of $B^2=[x_{i,j}^{(2)}]$ we have the strict inequalities $\sum_{i=1}^{d}x_{i,j}^{(2)}<1$ for all $j$.
\end{lemma}
\begin{proof}
Set $\epsilon_j:=\sum_{i=1}^{d}b_{i,j}$ for all $j\geq1$.
Then $\epsilon_j\in(0,1]$, and we fix an index $j_0$ such that $\epsilon_{j_0}<1$.
Now set $\tilde{b}_{i,j}:=\epsilon_j^{-1}b_{i,j}$ and $B:=[b_{i,j}]$;
then $B=\tilde{B}D$, where $\tilde{B}$ has all column sums one and $D$ is a diagonal matrix with $(j,j)$-entry $\epsilon_j$. Then, since $B$ is irreducible, we obtain:
\begin{equation}\label{A.1}
x_{i,j}^{(2)}=\sum_{k=1}^{d}\tilde{b}_{i,k}\tilde{b}_{k,j}\epsilon_k\epsilon_j\leq\sum_{k=1}^{d}\tilde{b}_{i,k}\tilde{b}_{k,j}\epsilon_k<\sum_{k=1}^{d}\tilde{b}_{i,k}\tilde{b}_{k,j}=\tilde{x}_{i,j}^{(2)}.
\end{equation}
given that $x_{i,j}^{(2)}>0$. (Otherwise $x_{i,j}^{(2)}=\tilde{x}_{i,j}^{(2)}=0$.)
Summing over $i$, we arrive at the required conclusion.
\end{proof}

\begin{corollary}\label{lemA}
Let $B\in\mathbb{R}^{d\times d}$ be an irreducible matrix with nonnegative real entries $[b_{i,j}]$ satisfying the system of inequalities $\sum_{i=1}^{d}b_{i,j}\leq1$. If at least one of these inequalities is strict, then the sequence $B^n$ converges to the zero matrix.
\end{corollary}

\begin{proof}
Indeed, let $\alpha=\max\limits_{i,j=1,2,\dots,d}\big\{x_{i,j}^{(2)}/\tilde{x}_{i,j}^{(2)}\big\}$. (If $x_{i,j}^{(2)}=\tilde{x}_{i,j}^{(2)}=0$, then we set $0/0=0$.) It follows from \eqref{A.1} that $\alpha<1$. Then $B^2\leq\alpha\tilde{B}^2$ entrywise. Consequently, $B^{2n}\leq\alpha^n\tilde{B}^{2n}$, and as $n\to\infty$, $B^{2n}$ converges to the zero matrix. As well, $B^{2n+1}=B^{2n}B$ must converge to the zero matrix as $n\to\infty$.
\end{proof}

\section*{ACKNOWLEDGEMENTS.} The author thanks all the people, who made critical comments officially or privately.

\subsection*{Disclosures and declarations.}
No conflict of interests was reported by the author.

\subsection*{Declaration of funding.}
No funding for this research was received.

\section*{ORCID}
Vyacheslav M. Abramov http://orcid.org/0000-0002-9859-100X

\bibliographystyle{vancouver}

\vfill\eject

\end{document}